\documentclass{amsart}
\usepackage{amssymb,latexsym}

\newcommand{\doublesubscript}[3]{
\displaystyle\mathop{\displaystyle #1_{#2}}_{#3}
}

\numberwithin{equation}{section}

\newtheorem{theorem}{Theorem}[section]
\newtheorem{proposition}[theorem]{Proposition}
\newtheorem{corollary}[theorem]{Corollary}
\newtheorem{definition}[theorem]{Definition}

\newtheorem{remark}[theorem]{Remark}
\newtheorem{lemma}[theorem]{Lemma}

\newtheorem{example}[theorem]{Example}

\def\mathbb\bf
\def\mathbf\bf
\def\mathcal\cal

\def\ZZ{{\bf Z}}

\def\FF{{\mathcal F}}
\def\II{{\mathcal I}}
\def\CCC{{\mathcal C}}
\def\Fl{\mbox{\rm Fl}}

\def\leqdeg{{\stackrel{deg}{\leq}}}

\def\aa{{\mathbf{a}}}
\def\bb{{\mathbf{b}}}
\def\cc{{\mathbf{c}}}
\def\dd{{\mathbf{d}}}

\def\Del{\Delta}

\def\Lam{\Lambda}

\begin{document}

\title{Multiple flag varieties of finite type}

\author{Peter Magyar}
\address{Dept. of Mathematics,
Northeastern University\\
Boston, Massachusetts 02115, USA
\\ \mbox{}}
\email{pmagyar@lynx.neu.edu}

\author{Jerzy Weyman}
\email{weyman@neu.edu}

\author{Andrei Zelevinsky}
\email{andrei@neu.edu}

\thanks{The research of Peter Magyar,
Jerzy Weyman and Andrei Zelevinsky
was supported in part by an NSF Postdoctoral Fellowship and
NSF grants \#DMS-9700884 and
\#DMS-9625511, respectively.}


\date{\today }


\begin{abstract}

We classify all products of flag varieties with finitely many orbits
under the diagonal action of the general linear group.
We also classify the orbits in each case and construct
explicit representatives.

\end{abstract}

\maketitle

\section{Introduction}
\label{sec:intro}

For a reductive group $G$, the Schubert (or Bruhat) decomposition
describes the orbits of a Borel subgroup $B$ acting on
the flag variety $G/B$.  It is the starting point for analyzing
the geometry and topology of $G/B$, and is also significant for the
representation theory of $G$.
This decomposition states that $G/B = \coprod_{w \in W} B\cdot w B$,
where $W$ is the Weyl group.
An equivalent form which is more symmetric is
the $G$-orbit decomposition of the double flag variety:
$(G/B)^2 = \coprod_{w \in W} G\cdot(eB, wB)$.
One can easily generalize the Schubert decomposition by considering
$G$-orbits on a product of two {\it partial} flag varieties $G/P \times G/Q$,
where $P$ and $Q$ are parabolic subgroups.
The crucial feature in each case is that the number of orbits is {\it finite}
and has a rich combinatorial structure.

Here we address the more general question:
for which tuples of parabolic subgroups
$(P_1, \ldots, P_k)$ does the group $G$ have finitely many orbits
when acting diagonally in the product of several flag varieties
$G/P_1 \times \cdots \times G/P_k$?
As before, this is equivalent to asking when
$G/P_2\times\cdots\times G/P_{k}$ has finitely many $P_1$-orbits.
(If $P_1=B$ is a Borel subgroup, this is one definition
of a spherical variety.  Thus, our problem includes that
of classifying the multiple flag varieties of spherical type.)
To the best of our knowledge, the problem of classifying all
finite-orbit tuples for an arbitrary $G$ is still open
(although in the special case when $k = 3$, $P_1 = B$,
and $P_2$ and $P_3$ are maximal parabolic subgroups,
such a classification was given in \cite{littelmann}).

In this paper, we present a complete solution
of the classification problem for $G = GL_n$.
For this case, the partial flag varieties $G/P$ consist of all flags
of subspaces with some fixed dimensions in an $n$-dimensional
vector space $V$.
Our classification theorem (Theorem~\ref{th:finite-type-list})
provides the list of all dimension types such that
$GL_n$ has finitely many orbits
in the corresponding product of flag varieties.
We also classify the orbits in each case and construct
explicit representatives (standard forms).
Precise formulations of the main results will be given in the next
section; the proofs are given in Sections~\ref{sec:proofs-classification}
and \ref{sec:proofs-representatives}.
In Section~\ref{sec:Bruhat order} we also discuss some partial results
on the generalized Bruhat order given by adjacency of orbits.

We use results and ideas from the theory of quiver representations.
In fact, our key criterion for finite type (Proposition~\ref{pr:criterion}
below) is very close to (but distinct from) the characterization of
quiver representations of finite type due to V.~Kac~\cite{kac}.

\textsc{Acknowledgments.}
We thank Michel Brion and Claus Ringel for helpful
references and discussions.
This paper was completed during Andrei Zelevinsky's stay at
the Institut Fourier, St Martin d'H\`eres (France) in February - March 1998;
he gratefully acknowledges the warm hospitality of Michel Brion and
financial support by CNRS, France.

\section{Main results}
\label{sec:main}

\subsection{Classification theorem}
\label{sec:classification}

Let $\aa = (a_1, \ldots, a_p)$ be a nonnegative list of integers
with sum equal to $n$.
We call such a list a \emph{composition} of $n$,
and $a_1, \ldots, a_p$ the {\it parts} of $\aa$.
Thoughout this paper, all vector spaces are
over a fixed algebraically closed field.
For a vector space $V$ of dimension $n$, we denote by
${\rm Fl}_{\aa} (V)$ the variety of flags
$A = (0 = A_0 \subset A_1 \subset \cdots \subset A_p = V)$
of vector subspaces in $V$ such that
$$\dim (A_{i}/A_{i-1}) = a_i \quad (i = 1, \ldots, p) \ .$$
A tuple of compositions $(\aa_1, \ldots, \aa_k)$ of
the same number $n$ is said to be of \emph{finite type}
if the group $GL (V)$ (acting diagonally) has finitely many orbits in
the multiple flag variety
${\rm Fl}_{\aa_1} (V) \times \cdots \times {\rm Fl}_{\aa_k} (V)$.
We say that a composition is \emph{trivial}
if it has only one non-zero part $n$.
Then the corresponding flag variety consists of a single point,
so adding any number of trivial compositions to a tuple
gives essentially the same multiple flag variety, and
does not affect the finite type property.

\begin{theorem}
\label{th:k=3}
If a tuple of non-trivial compositions $(\aa_1, \ldots, \aa_k)$ is
of finite type then $k \leq 3$.
\end{theorem}

In other words, a multiple flag variety of finite type cannot have more
than $3$ non-trivial factors.
Thus any tuple of compositions of finite type can be
made into a triple by adding or removing trivial compositions,
and we need only classify \emph{triples} of finite type.
We will write $(\aa, \bb, \cc)$ instead of
$(\aa_1, \aa_2, \aa_3)$.
We denote by ${\rm min} (\aa)$ the minimum of the non-zero
parts of a composition $\aa$.
Now we can formulate our first main theorem.

\begin{theorem}
\label{th:finite-type-list}
Let $\aa, \bb$ and $\cc$ be compositions and
let $p,q$ and $r$ denote their respective numbers of {\em non-zero} parts.
Assume without loss of generality that $p \leq q \leq r$.
Then the triple $(\aa, \bb, \cc)$ is of finite type
if and only if it belongs to one of the following classes:

\begin{tabbing}
{\rm ($D_{r+2}$)\ }\=$(p,q,r) = (2,2,r)$, $2 \leq r$.
\kill
{\rm ($A_{q,r}$)\ }\> $(p,q,r) = (1,q,r)$, $1 \leq q \leq r$.
\\[.5em]
{\rm ($D_{r+2}$)\ }\>$(p,q,r) = (2,2,r)$, $2 \leq r$.
\\[.5em]
{\rm ($E_6$)\  }\>$(p,q,r) = (2,3,3)$.
\\[.5em]
{\rm ($E_7$)\ }\> $(p,q,r) = (2,3,4)$.
\\[.5em]
{\rm ($E_8$)\ }\> $(p,q,r) = (2,3,5)$.
\\[.5em]
{\rm ($E^{(a)}_{r+3}$)\ }\>$(p,q,r) = (2,3,r)$, $3 \leq r$,
${\rm min} \ (\aa) = 2$.
\\[.5em]
{\rm ($E^{(b)}_{r+3}$)\ }\>$(p,q,r) = (2,3,r)$, $3 \leq r$,
${\rm min} (\bb) = 1$.
\\[.5em]
{\rm ($S_{q,r}$)\ }\>$(p,q,r) = (2,q,r)$, $2 \leq q \leq r$,
${\rm min} (\aa) = 1$.
\end{tabbing}
\end{theorem}

The types in Theorem~\ref{th:finite-type-list}
have some obvious overlaps.
The type $A_{q,r}$ covers all multiple flag varieties
with less than three non-trivial factors.
The type $S_{q,r}$ appeared in Brion \cite{brion}.
For relations with the classification
of quiver representations of finite type due to V.~Kac \cite{kac},
see Remark \ref{re:quiver classification}.

Note that, for each of the first five types in
Theorem~\ref{th:finite-type-list}, there are no restrictions
on the dimensions of subspaces
in the corresponding flag varieties; only the last three types
$E^{(a)}$, $E^{(b)}$ and $S$ involve such restrictions.
The first five types are naturally related to the simply-laced
Dynkin graphs (as suggested by their names).
Let $T=T_{p,q,r}$ denote the graph with
$p+q+r-2$ vertices that consists of $3$ chains with
$p,q$, and $r$ vertices, joined together at a common endpoint.
We see that the cases in our classification with no restrictions
on dimensions are precisely those for which $T$ is one of the
Dynkin graphs $A_n$, $D_n$, $E_6$, $E_7$, or $E_8$.
This is of course no coincidence: we will see that this part
of our classification is equivalent to Gabriel's classification of
quivers of finite type (and follows from the Cartan-Killing
classification of graphs that give rise to positive-definite quadratic forms).

\subsection{Classification of orbits}

Now we describe a combinatorial parametrization of the set of
$GL (V)$-orbits in
${\rm Fl}_{\aa} (V) \times {\rm Fl}_{\bb} (V)  \times {\rm Fl}_{\cc} (V)$
for any triple $(\aa, \bb, \cc)$ of compositions of finite type.
For a composition $\aa = (a_1, \ldots, a_p)$, we write
$$
|\aa| = a_1 + \cdots + a_p \quad ,\quad
\Vert \aa \Vert^2 = a_1^2 + \cdots + a_p^2 \ ;
$$
the number $p$ of parts of $\aa$ will be denoted $\ell (\aa)$
and called the \emph{length} of $\aa$.

For any positive integers $p, q$, and $r$, let
$\Lam_{p,q,r}$ denote the additive semigroup of all triples
of compositions $(\aa, \bb, \cc)$ such that
$(\ell (\aa), \ell (\bb), \ell (\cc)) = (p,q,r)$,
and $|\aa| = |\bb| = |\cc|$.
(Here, in contrast to the notation of Theorem~\ref{th:finite-type-list},
the numbers $p,q,r$ include the zero parts of $\aa$, $\bb$, $\cc$.)
For every $(\aa, \bb, \cc) \in \Lam_{p,q,r}$, we set
$$
Q(\aa, \bb, \cc) = \dim GL (V) - \dim {\rm Fl}_{\aa} (V)
- \dim {\rm Fl}_{\bb} (V) - \dim {\rm Fl}_{\cc} (V) \ ,
$$
where $V$ is a vector space of dimension $n = |\aa| = |\bb| = |\cc|$.
An easy calculation shows that
\begin{equation}
\label{eq:Q-formula}
Q(\aa, \bb, \cc) = \frac{1}{2}
(\Vert \aa \Vert^2 + \Vert \bb \Vert^2 + \Vert \cc \Vert^2 - n^2) \ .
\end{equation}
The function $Q$ is called the \emph{Tits quadratic form}.

Let $\Pi_{p,q,r}$ denote the set of all triples
$\dd = (\aa, \bb, \cc) \in \Lam_{p,q,r}$ of finite type such that
$Q(\dd) = 1$.

\begin{theorem}
\label{th:orbits}
Let $(\aa, \bb, \cc) \in \Lam_{p,q,r}$ be a triple of compositions
of finite type.
Then $GL (V)$-orbits in
${\rm Fl}_{\aa} (V) \times {\rm Fl}_{\bb} (V)  \times {\rm Fl}_{\cc} (V)$
are in natural bijection with families
of nonnegative integers $M = (m_{\dd})$
indexed by $\dd \in \Pi_{p,q,r}$ such that,
in the semigroup $\Lam_{p,q,r}$,
$$\sum_{\dd \in \Pi_{p,q,r}} m_{\dd} \dd = (\aa, \bb, \cc).$$
\end{theorem}

The set $\Pi_{p,q,r}$  can be explicitly described as follows.
For a composition $\aa$, we denote by $\aa^+$ the partition
obtained from $\aa$ by removing all zero parts and
rearranging the non-zero parts in weakly decreasing order.
(For example, if $\aa = (0, 2, 1, 0, 3, 2)$ then $\aa^+ = (3,2,2,1)$.)
We denote
$(a^p)=\underbrace{(a,\ldots,a)}_{p\ \mbox{\tiny parts}}$.

\begin{theorem}
\label{th:primitive}
A triple $(\aa, \bb, \cc) \in \Lam_{p,q,r}$ belongs to $\Pi_{p,q,r}$
if and only if the (unordered) triple of partitions
$\{\aa^+, \bb^+, \cc^+\}$ is one of the following:
$$
\{(1),(1),(1)\},\quad \{(3^2),(2^3),(2,1,1,1,1)\},
\quad \{(4,2),(2^3),(1^6)\},
$$
$$
\{(m+1,m),(m,m,1),(1^{2m+1})\} \ \ (m \geq 2),\quad
\{(m,m),(m,m-1,1),(1^{2m})\} \ \ (m \geq 2),
$$
$$
\{(n-1,1),(1^n),(1^n)\} \  (n \geq 2).
$$
\end{theorem}

\begin{remark}
\label{re:Simpson}
{\rm Except for $((3^2),(2^3),(2,1,1,1,1))$,
all of the triples in Theorem~\ref{th:primitive} are ``spherical'', meaning
that one of the compositions $\aa,\bb,\cc$ is equal to $(1^n)$,
or equivalently one of the factors of the triple flag variety
is a complete flag variety.
Our spherical cases are identical to C. Simpson's
list \cite{simpson} of certain local systems on ${\bf P}^1$
with three punctures.
Such local systems are equivalent to triples of
matrices $X_1,X_2,X_3 \in GL_n$ with $X_1 X_2 X_3=I$,
up to simultaneous conjugation.
Simpson classifies the triples of semi-simple conjugacy classes
$C_1, C_2, C_3$ of $GL_n$, one of which is regular,
such that the product $C_1\times C_2\times C_3$
contains a unique solution to the equation $X_1 X_2 X_3 = I$
(up to simultaneous conjugation).  These are called {\it rigid}
local systems.  The last indecomposable type on our list corresponds
to the local system associated to the Pochhammer hypergeometric function.
}
\end{remark}
 
To describe the bijection in Theorem~\ref{th:orbits}, we
introduce the following additive category $\FF_{p,q,r}$.
The \emph{objects} of ${\FF}_{p,q,r}$ are families $(V; A, B, C)$,
where $V$ is a finite-dimensional vector space, and
$(A,B,C)$ is a triple of flags in $V$ belonging to
${\rm Fl}_{\aa} (V) \times {\rm Fl}_{\bb} (V)  \times {\rm Fl}_{\cc} (V)$
for some $(\aa, \bb, \cc) \in \Lam_{p,q,r}$.
The triple $\dd = (\aa, \bb, \cc)$ is called the \emph{dimension vector} of
$(V; A, B, C)$.
A \emph{morphism} from $(V; A, B, C)$ to $(V'; A', B', C')$
in $\FF_{p,q,r}$ is a linear map $f: V \to V'$ such that
$f(A_i) \subset A'_i, \, f(B_i) \subset B'_i$, and $f(C_i) \subset C'_i$
for all $i$.
Direct sum of objects is taken componentwise
on each member of each flag in the objects.

Comparing definitions, we see that isomorphism classes of objects
in $\FF_{p,q,r}$ with a given dimension vector
$(\aa, \bb, \cc) \in \Lam_{p,q,r}$ are naturally
identified with $GL (V)$-orbits in
${\rm Fl}_{\aa} (V) \times {\rm Fl}_{\bb} (V)  \times {\rm Fl}_{\cc} (V)$.
The advantage of dealing with $\FF_{p,q,r}$ is that this category
admits direct sums, and so every object $(V; A, B, C)$ of $\FF_{p,q,r}$ can be
decomposed into a direct sum of indecomposable objects.
By the Krull-Schmidt theorem (see \S\ref{sec:indecomposables}),
such a decomposition is unique up to an automorphism of $(V; A, B, C)$.
So the isomorphism class of an object is determined by
the multiplicities of the non-isomorphic indecomposable objects
in its decomposition.
Theorem~\ref{th:orbits} now becomes a consequence of the following.

\begin{theorem}
\label{th:indecomposables}
For every $\dd \in \Pi_{p,q,r}$, there exists a unique
isomorphism class $I_{\dd}$ of indecomposable objects in $\FF_{p,q,r}$
with the dimension vector $\dd$.
For every triple $(\aa, \bb, \cc) \in \Lam_{p,q,r}$
of finite type, any object in $\FF_{p,q,r}$ with the dimension
vector $(\aa, \bb, \cc)$ decomposes (uniquely) into a direct sum of objects
$I_{\dd}$.
\end{theorem}

\begin{corollary}
\label{cor:orbit parametrization}
The bijection in Theorem~\ref{th:orbits}
sends a family $M = (m_{\dd}) \,\, (\dd \in \Pi_{p,q,r})$ to the
$GL(V)$-orbit $\Omega_M$ in
${\rm Fl}_{\aa} (V) \times {\rm Fl}_{\bb} (V)  \times {\rm Fl}_{\cc} (V)$
corresponding to the isomorphism class
$\bigoplus_{\dd \in \Pi_{p,q,r}} m_{\dd} I_{\dd}$ of objects in
$\FF_{p,q,r}$.
\end{corollary}

\subsection{Representatives of orbits}

By Corollary~\ref{cor:orbit parametrization},
in order to give an explicit representative of each $GL(V)$-orbit
in a multiple flag variety of finite type, it is enough to exhibit
a triple of flags that represents every indecomposable object
$I_{\dd}$ in Theorem~\ref{th:indecomposables}.
All possible dimension vectors $\dd = (\aa, \bb, \cc) \in \Pi_{p,q,r}$
are described in Theorem~\ref{th:primitive}.
Note that vanishing of some part $a_i$ in a composition $\aa$
means that in any flag $A \in {\rm Fl}_{\aa} (V)$ the subspace
$A_i$ coincides with $A_{i-1}$.
Thus in constructing a representative for $I_{\dd}$, we can assume without
loss of generality that none of the compositions $\aa, \bb$ and $\cc$
have zero parts.
So it is enough to treat all the dimension vectors $\dd$ obtained from
the triples $(\aa^+, \bb^+, \cc^+)$ in Theorem~\ref{th:primitive}
by some permutations of the parts.
In particular, we assume that $p = \ell (\aa) \leq 2$;
thus a flag $A \in {\rm Fl}_{\aa} (V)$ is determined by one vector
subspace $A_1$
in $V$.
Under these assumptions, we will show that a triple of flags $(V;A,B,C)$
representing $I_\dd$ can be presented in a standard form
according to the following definition.

\begin{definition}
\label{def:standard form}
{\rm An object $(V; A, B, C)$ in ${\FF}_{p,q,r}$
is presented in {\it standard form} if
$V$ is given a basis $e_1, \ldots, e_n$ with the
following properties:

\smallskip

\noindent (1) Each subspace $B_i$ of the flag $B$ has a basis
consisting of the first
$b_1 + \cdots + b_i$ standard basis vectors $e_1, e_2, \ldots$,
while each $C_i$ has a basis consisting of the last $c_1 + \cdots + c_i$
basis vectors $e_n, e_{n-1}, \ldots$.

\smallskip

\noindent (2) $p \leq 2$, and the vector subspace
$A_1 \subset V$ of dimension $a = a_1$ has basis vectors
$\sum_{l \in S_1} e_l, \ldots$, $\sum_{l \in S_{a}} e_l$
for some subsets $S_1, \ldots, S_{a} \subset \{1, \ldots, n\}$.

\smallskip

\noindent (3) The subsets $S_k$ satisfy:
$|\bigcup_{k \neq k'} (S_k \cap S_{k'})| \leq 2 \ .$

}
\end{definition}

\begin{theorem}
\label{th:standard form}
Let $\dd \in \Pi_{p,q,r}$ be a triple of compositions
obtained from some triple $(\aa^+, \bb^+, \cc^+)$ in
Theorem~\ref{th:primitive} by permutations of the parts.
Then the corresponding indecomposable object $I_{\dd}$ in ${\FF}_{p,q,r}$
can be presented in standard form with the
collection of subsets $S_1, \ldots, S_{a}$ chosen as follows:
\\[.4em]
$\dd = ((1),(1),(1)):
\quad S_1= \{1\}.$
\\[.4em]
$\dd =((4,2),(2^3),(1^6)):
\quad S_1= \{1,5\}, \  S_2 = \{2,3,\}, \ S_3 = \{2,5,6\}, \  S_4= \{4,5\}.$
\\[.2em]
$\dd =((3^2),(2^3),(2,1,1,1,1)):
\quad S_1= \{1,2,3\}, \  S_2= \{1,6\}, \  S_3= \{2,4,5\}.$
\\[.2em]
$\dd =((3^2),(2^3),(1,1,2,1,1)):
\quad S_1= \{1,5,6\}, \  S_2= \{2,3,6\}, \  S_3= \{4,5\}.$
\\[.2em]
$\dd =((3^2),(2^3),(1,1,1,1,2)):
\quad S_1= \{1,4,6\},\   S_2= \{2,4,5\}, \  S_3= \{2,3\}.$
\\[.4em]
$\dd =((3^2),(2^3),(1,2,1,1,1)):
\quad S_1= \{1,2,4,6\}, \  S_2= \{1,3\}, \  S_3= \{1,5\}.$
\\[.2em]
$\dd =((3^2),(2^3),(1,1,1,2,1)):
\quad S_1= \{1,2,4,6\}, \  S_2= \{1,3\}, \  S_3= \{1,5\}.$
\\[.2em]
$\dd =((2,4),(2^3),(1^6)):
\quad S_1 = \{1,2,3,6\} , \  S_2 = \{1,4,5\}.$
\\[.4em]
$\dd = ((m,m+1),(1,m,m),(1^{2m+1})):
\\
\mbox{}
\quad S_k = \{1,k+1,2m+2-k\} \ \  (1 \leq k  \leq m). $
\\[.2em]
$\dd = ((m+1,m),(1,m,m),(1^{2m+1})):
\\
\mbox{}
\quad S_1 = \{1,2\}, \  S_k = \{1,k+1,2m+3-k\} \ \  (2 \leq k \leq m),
\quad S_{m+1} = \{1,m+2\}. $
\\[.2em]
$\dd = ((m,m+1),(m,1,m),(1^{2m+1})):
\\
\mbox{}
\quad S_k = \{k,m+1,2m+2-k\} \ \  (1 \leq k \leq m).$
\\[.2em]
$\dd = ((m+1,m),(m,1,m),(1^{2m+1})):
\\
\mbox{}
\quad S_{1} = \{1,m+1\}, \  S_k = \{k,m+1,2m+3-k\} \ \  (2 \leq k \leq m),
\  S_{m+1} = \{m+1,m+2\}. $
\\[.2em]
$\dd = ((m,m+1),(m,m,1),(1^{2m+1})):
\\
\mbox{}
\quad S_k = \{k,2m+1-k,2m+1\} \ \  (1 \leq k \leq m). $
\\[.2em]
$\dd = ((m+1,m),(m,m,1),(1^{2m+1})): \quad S_1 = \{1,2m+1\},
\\
\mbox{} \quad S_k = \{k,2m+2-k,2m+1\} \ \  (2 \leq k \leq m),
\  S_{m+1} = \{m+1, 2m+1\}. $
\\[.4em]
$\dd = ((m,m),(1,m-1,m),(1^{2m})):
\\
\mbox{}
\quad S_k = \{1,k+1,2m+1-k\} \ \  (1 \leq k \leq m-1),
\  S_{m} = \{1,m+1\}. $
\\[.2em]
$\dd = ((m,m),(1,m,m-1),(1^{2m})):
\\
\mbox{}
\quad S_{1} = \{1,2\}, \  S_k = \{1,k+1,2m+2-k\} \ \
(2 \leq k \leq m). $
\\[.2em]
$\dd = ((m,m),(m-1,1,m),(1^{2m})):
\\
\mbox{}
\quad S_k = \{k,m,2m+1-k\} \ \
(1 \leq k \leq m-1), \  S_{m} = \{m,m+1\}. $
\\[.2em]
$\dd = ((m,m),(m,1,m-1),(1^{2m})):
\\
\mbox{}
\quad S_{1} = \{1,m+1\}, \  S_k = \{k,m+1,2m+2-k\} \ \
(2 \leq k \leq m). $
\\[.2em]
$\dd = ((m,m),(m-1,m,1),(1^{2m})):
\\
\mbox{}
\quad S_k = \{k,2m-k,2m\} \ \
(1 \leq k \leq m-1), \  S_{m} = \{m,2m\}. $
\\[.2em]
$\dd = ((m,m),(m,m-1,1),(1^{2m})):
\\
\mbox{}
\quad S_{1} = \{1,2m\}, \  S_k = \{k,2m+1-k,2m\} \ \
(2 \leq k \leq m). $
\\[.4em]
$\dd = ((n-1,1),(1^n),(1^n)):
\quad S_k = \{1,k+1\} \ \
(1 \leq k \leq n-1).$
\\[.2em]
$\dd = ((1,n-1),(1^n),(1^n)):
\quad S_1 = \{1,2, \ldots, n\}.$
\end{theorem}

\vspace{1em}

For any composition $\aa$,
let $\aa_{\rm red}$ denote the composition obtained from $\aa$
by removing all zero parts and
keeping the non-zero parts in the same order.
For a dimension vector $\dd = (\aa, \bb, \cc) \in \Lam_{p,q,r}$, we set
$\dd_{\rm red} = (\aa_{\rm red},\bb_{\rm red},\cc_{\rm red})$
and call $\dd_{\rm red}$ the \emph{reduced dimension vector} of $\dd$.
Thus $\dd \in \Pi_{p,q,r}$ if and only if $\dd_{\rm red}$
is one of the triples in Theorem~\ref{th:standard form}.

Finite-type indecomposables in Theorem~\ref{th:standard form}
are closely related to \emph{exceptional pairs} in the sense 
of \cite{ringel}; see Remark \ref{re:mutations}.

\begin{example}
\label{ex:A-orbits}
{\rm
Type $A_{q,r}$: two flags.
Let $\bb = (b_1, \ldots, b_q)$ and $\cc = (c_1, \ldots, c_r)$
be two compositions of $n$.
We can identify a pair of partial flags
$(B,C) \in {\rm Fl}_{\bb} (V)  \times {\rm Fl}_{\cc} (V)$
with the object $(V;A,B,C)$ in the category $\FF_{1,q,r}$,
where $A$ is the trivial flag $(0 = A_0 \subset A_1 = V)$.
An indecomposable summand of $(V;A,B,C)$ can only
have the reduced dimension vector $((1),(1),(1))$.
The indecomposable objects with this reduced dimension vector are of the form
$I_{ij} = (V';A',B',C')$ where $1\leq i \leq q$ and $1\leq j\leq r$:
here $\dim V'=1$,\  $A'$ is the trivial flag in $V'$,
$B' = (0 =B'_0 = \cdots = B'_{i-1} \subset B'_i = \cdots = B'_q = V')$,
and $C' = (0=C'_0 = \cdots = C'_{j-1} \subset C'_j = \cdots = C'_r = V')$.

It follows that $GL(V)$-orbits in
${\rm Fl}_{\bb} (V)  \times {\rm Fl}_{\cc} (V)$
are parametrized by $q \times r$ nonnegative
integer matrices $M = (m_{ij})$ with row sums
$b_1, \ldots, b_q$ and column sums $c_1, \ldots, c_r$;
the orbit $\Omega_M$ corresponds to a direct sum
$\bigoplus_{i,j} m_{ij} I_{ij}$.
(In particular, if $B$ and $C$ are complete flags,
we obtain the usual parametrization of orbits by permutation matrices.)
A representative of $\Omega_M$ can be given as follows:
$V$ has a basis
$\{e_{ijk} :  1\leq i \leq q, \  1\leq j\leq r, \
1\leq k\leq m_{ij}\}$,
each $B_{i}$ is spanned by the $e_{i'j'k'}$ with $i' \leq i$,
and each $C_{j}$ is spanned by the $e_{i'j'k'}$ with $j' \leq j$.
}
\end{example}

\begin{example}
\label{ex:S-orbits}
{\rm Type $S_{q,r}$: two flags and a line.
As in Example~\ref{ex:A-orbits}, let $\bb$ and $\cc$ be any two compositions
of $n$, but now let us take $\aa = (1,n-1)$.
Let $(V;A,B,C)$ be a triple of flags of type $(S_{q,r})$
with the dimension vector $(\aa, \bb, \cc) \in \Lam_{2,q,r}$.
By inspection of the cases in Theorem~\ref{th:standard form},
we see that an indecomposable summand of $(V;A,B,C)$ can only
have the reduced dimension vector $((1),(1),(1))$ or
$((1,t-1),(1^t),(1^t))$ for some $t = 2\  \ldots, n$.
The corresponding indecomposable objects are of the following form.
First each $I_{ij}$ in the previous example can be also
considered as an indecomposable object in the present situation:
we take $V', B'$, and $C'$ as above, and define the flag
$A'$ as $(0 = A'_0 = A'_1 \subset A'_2 = V')$;
by abuse of notation, we denote this indecomposable object
in $\FF_{2,q,r}$ by the same symbol $I_{ij}$.

Besides these indecomposables, the object $(V;A,B,C)$
must have precisely one indecomposable summand $(V';A',B',C')$
with $\dim A'_1 = 1$ (since $\dim A_1 = 1$).
Such indecomposables are indexed by non-empty sets
$\Del = \{ (i_1,j_1), (i_2,j_2),$ $\ldots ,$ $(i_{t}, j_{t})\}$
with $1\leq i_1<\ldots <i_t \leq q$ and
$r \geq j_1 > \ldots > j_t \geq 1$
(such a $\Del$ can be pictured as the outer corners of
a Young diagram contained in a $q \times r$ rectangle).
The indecomposable object $I_{\Del}$ is represented by the following
triple of flags $(V';A',B',C')$: the space $V'$ has basis
$e_1,\ldots,e_t$,\ the subspace $A'_1$ is spanned by
$e_1+\ldots+e_t$,\ each $B'_{i}$ is spanned by the $e_l$
with $i_l \leq i$,\ and each $C'_{j}$ is spanned by the $e_l$
with $j_l \leq j$.

We see that $GL(V)$-orbits in
${\rm Fl}_{\aa} (V) \times {\rm Fl}_{\bb} (V)  \times {\rm Fl}_{\cc} (V)$
correspond to objects $I_{\Del} \oplus {\bigoplus} m'_{ij} I_{ij}$
in $\FF_{2,q,r}$  with the dimension vector $(\aa, \bb, \cc)$.
It is convenient to use the numbers
$m_{ij}=m'_{ij}+1$ for $(i,j)\in \Del$ and
$m_{ij}=m'_{ij}$ otherwise.
In this notation, the orbits are parametrized by pairs $(\Del, M = (m_{ij}))$
where $\Del$ is any set as above, and $M$ is a $q \times r$
nonnegative integer matrix with row sums
$b_1, \ldots, b_q$ and column sums $c_1, \ldots, c_r$,
and with $m_{ij}>0$ for all $(i,j)\in \Del$.
A representative of the orbit $\Omega_{\Del,M}$ is
$(V; A,B,C)$ where $(V;B,C)$ is the representative of
$\Omega_M$ constructed in the previous example, and
$A_1$ is spanned by the vector $\sum_{(i,j) \in \Del} e_{ij1}$.
}
\end{example}

\section{Proofs of Theorems~\ref{th:k=3} --- \ref{th:indecomposables}}
\label{sec:proofs-classification}

For a $k$-tuple of positive integers $(p_1,\ldots,p_k)$,
we define the graph $T_{p_1,\ldots,p_k}$,
the semigroup $\Lam_{p_1,\ldots,p_k}$, the Tits quadratic form
$Q(\dd)$ on $\Lam_{p_1,\ldots,p_k}$,
and the additive category $\FF_{p_1,\ldots,p_k}$
analogously to their counterparts for $k=3$.
Also let $\Pi_{p_1,\ldots,p_k}$ be the set of $k$-tuples of compositions
of finite type $\dd \in \Lam_{p_1,\ldots,p_k}$ with $Q(\dd)=1$.
When there is no risk of ambiguity, we drop the subscripts
$(p_1,\ldots,p_k)$ and write $\Lam$, $\FF$, $\Pi$, etc.

\subsection{Proof of Theorem~\ref{th:indecomposables}}
\label{sec:indecomposables}

The only ``non-elementary" part of our argument is the following
proposition.

\begin{proposition}
\label{pr:Kac indecomposables}
Suppose $\dd \in \Lam$ is the dimension vector of an indecomposable
object of $\FF$ with $Q(\dd) \geq 1$.
Then $Q(\dd) = 1$, and there is a unique
isomorphism class $I_{\dd}$ of indecomposable objects
with the dimension vector $\dd$.
\end{proposition}

\noindent {\it Proof.} Consider $T=T_{p_1,\ldots,p_k}$ as a
directed graph with all edges pointing toward the central vertex
where all the chains are joined.
Let $\CCC$ be the category of quiver representations of $T$:
recall that such a representation is specified by attaching
a finite-dimensional vector space to every vertex of $T$, and
a linear map between the corresponding spaces to every arrow
(directed edge) of $T$.
There is an obvious functor from $\FF$ to $\CCC$: given a tuple of
flags, the corresponding quiver representation associates
the flag subspaces to vertices of $T$, and inclusion maps to arrows.
The image of this functor lies in the subcategory $\II$ of $\CCC$
consisting of quiver representations with all arrows represented by
injective linear maps.
In fact, our functor allows us to identify the isomorphism
classes of indecomposable objects in $\FF$
with those in $\II$.  Note that $\II$ is a full additive
subcategory of $\CCC$, and that
indecomposables of $\II$ are also indecomposables of $\CCC$,
since an injective linear map can never have a
non-injective map as a direct summand.

In view of this translation, our proposition follows from
general results due to V.~Kac (\cite[Theorem 1]{kac}) which provide
a description of dimension vectors for indecomposable
quiver representations of an arbitrary finite directed graph.
(These dimension vectors turn out to be in a natural bijection
with positive roots of the simply-laced Kac-Moody Lie algebra
corresponding to the graph.)  Kac shows that the dimension vectors
of indecomposables all have $Q(\dd)\leq 1$; and if an indecomposable
has $Q(d)=1$ (the case of a real root), then there is a unique
indecomposable of dimension $\dd$ up to isomorphism.
This directly implies our Proposition.
%
%
%
%
%
%
%
%
\endproof \smallskip

Note that $\FF$ is not an abelian category (since it does
not always admit quotients). However, the Krull-Schmidt Theorem
(as in \cite{atiyah}) still applies.  That is,
each object of $\FF$ has a unique splitting into indecomposables.

In general, the condition that $\dd \in \Lam$ has $Q(\dd) = 1$
does not imply the existence of an indecomposable object $I_\dd$ in $\FF$
with the dimension vector $\dd$.
However, we now show that if $\dd$ is of finite type
with $Q(\dd)=1$ (i.e., if $\dd \in \Pi$)
then $I_\dd$ exists and has an important additional property.
For any two (isomorphism classes of) objects $F$ and $F'$ in $\FF$,
let us denote
\begin{equation}
\label{eq:dim Hom}
\langle F', F \rangle = \dim   {\rm Hom}_{\FF} (F', F) \ .
\end{equation}
We say that $F \in \FF$ is a \emph{Schur indecomposable} if
$\langle F , F \rangle = 1$ (which clearly implies that $F$ is indeed an
indecomposable object in $\FF$).

\begin{proposition}
\label{pr:finite indecomposables}
If $\dd \in \Pi$ then there exists a Schur indecomposable $I_\dd$
with the dimension vector $\dd$.
\end{proposition}

\noindent {\it Proof.}
Since $\dd$ is of finite type, the corresponding
multiple flag variety
$\Fl_\dd (V) = \Fl_{\aa_1}(V)  \times \cdots \times \Fl_{\aa_k} (V)$
has a (dense) Zariski open orbit $\Omega$.
Let $I_\dd$ be the corresponding isomorphism class in $\FF$, and let
$F$ be any representative of $I_\dd$; by abuse of notation, we
can think of $F$ as a point in $\Omega$.
Then we have
$$\langle I_\dd, I_\dd \rangle = \dim \mbox{Stab}_{GL(V)} (F) =
\dim GL(V) - \dim \Fl_\dd (V) = Q(\dd) = 1 \ .$$
Therefore, $I_\dd$ is a Schur indecomposable, as desired.
\endproof \smallskip

By Propositions~\ref{pr:Kac indecomposables} and
\ref{pr:finite indecomposables}, for every $\dd \in \Pi$,
there exists a unique isomorphism class $I_{\dd}$ of
indecomposable objects in $\FF$ with the dimension vector $\dd$.
Now the proof of Theorem~\ref{th:indecomposables}
(and hence that of Theorem~\ref{th:orbits} and
Corollary~\ref{cor:orbit parametrization}) can be concluded as
follows.

We say that a non-zero $\dd' \in \Lam$ is a \emph{summand} of
$\dd \in \Lam$ if $\dd - \dd' \in \Lam$.
It follows from the Krull-Schmidt theorem that if $\dd$ is of finite type
then every summand of $\dd$ is also of finite type.
Thus every object in $\FF$ whose dimension vector is of finite type
decomposes (uniquely) into a direct sum of objects $I_\dd$ for
$\dd \in \Pi$, and we are done.
\endproof \smallskip

\subsection{Proof of Theorems~\ref{th:k=3}, \ref{th:finite-type-list},
and \ref{th:primitive}}

The following criterion reduces the classification of tuples of compositions
of finite type to an ``elementary" problem about the Tits form.

\begin{proposition}
\label{pr:criterion}
A tuple of compositions $\dd \in \Lam$ is of
finite type if and only if $Q(\dd')\geq 1$ for any
summand $\dd'$ of $\dd$.
\end{proposition}

\noindent {\it Proof.} Let $\Fl_\dd (V)$ be the multiple flag variety
corresponding to $\dd \in \Lam$.
First suppose $\dd$ is of finite type.
Since the one-dimensional subgroup of scalar matrices in $GL(V)$
acts trivially on $\Fl_\dd (V)$, we must have
$\dim GL(V) - 1 \geq \dim \Fl_\dd (V)$, i.e., $Q(\dd)\geq 1$.
We have already noticed that any summand $\dd'$ of $\dd$
must also be of finite type, hence we must have $Q(\dd') \geq 1$.

Conversely, suppose $Q(\dd') \geq 1$ for any summand $\dd'$ of $\dd$.
In view of Proposition~\ref{pr:Kac indecomposables},
this implies that every indecomposable summand of an object
in $\FF$ with the dimension vector $\dd$ is uniquely determined
by its dimension vector.
Therefore, the isomorphism classes of objects with the dimension vector
$\dd$ are in a bijection with partitions of $\dd$ into
the sum of dimension vectors of indecomposables.
Since there are finitely many such partititons,
$\dd$ must be of finite type, and we are done.
\endproof \smallskip

\begin{remark}
\label{re:quiver classification}
{\rm  The criterion in Proposition \ref{pr:criterion} is almost identical
to that of V.~Kac \cite[Proposition 2.4]{kac}, for finite-type quiver
varieties: the quiver variety with a given dimension
vector $\dd$ has finitely many orbits exactly if $Q(\dd')\geq 1$
for all {\it quiver} summands $\dd'$ of $\dd$.
(A quiver summand need not have positive jumps in dimension along each flag;
only the dimension of each space must be positive.)

To illustrate the difference between the two criteria, consider
the triple of compositions $\dd=((3,1), (1^4), (1^4))$ of finite type
$S_{4,4}$.
The corresponding space of quiver representations
has infinitely many orbits because it has the quiver summand
$\dd' = ((2,2),(1^4),(1^4))$ with $Q(\dd') = 0$;
in terms of dimensions rather than dimension jumps,
{\footnotesize $\begin{array}{c@{\!}c@{\!}c}
123\,&4&\,321\\[-.4em] &3& \end{array}$ =
 $\begin{array}{c@{\!}c@{\!}c}
123\,&4&\,321\\[-.4em] &2& \end{array}$ +
 $\begin{array}{c@{\!}c@{\!}c}
000\,&0&\,000\\[-.4em] &1& \end{array}$
 }.
In fact, the infinitely many generic objects of dimension $\dd'$
are all subobjects of the unique generic
indecomposable $I_{\dd}$, but there are no quotient objects
in the flag category $\FF$, only in the quiver category $\CCC$.
(We thank C.~Ringel for this example.)
}
\end{remark}

The classification of finite types
in the rest of this section closely follows
the Cartan-Killing classification, which in our terminology amounts to
finding all graphs $T_{p_1, \ldots, p_k}$
with positive-definite Tits form.

\smallskip

\noindent {\it Proof of Theorem~\ref{th:k=3}.}
Clearly, if $\dd$ is a $k$-tuple of compositions of finite type
then any subtuple of $\dd$ is also of finite type.
Thus it suffices to show that a quadruple $\dd$ of non-trivial
compositions cannot be of finite type.
But any such quadruple has a summand
$\dd'$ with the reduced dimension vector $((1^2),(1^2),(1^2),(1^2))$.
A calculation shows that $Q(\dd') = 0$, so by
Proposition~\ref{pr:criterion}, $\dd$ cannot be of finite type.
\endproof \smallskip

Next, beginning the proof of Theorem~\ref{th:finite-type-list},
we eliminate those dimension vectors with a summand
corresponding to the minimal imaginary root of an
affine root system.
Let $N_{p,q,r}$ be the set of all
$\dd'= (\aa',\bb',\cc') \in \Lam_{p,q,r}$ such that
$\{\aa'_{\rm red},\bb'_{\rm red},\cc'_{\rm red}\}$
is one of the following three triples:
\begin{equation}
\label{eq:imaginary}
\{(1^3),(1^3),(1^3)\},\quad \{(2^2),(1^4),(1^4)\},\quad
\{(3^2),(2^3),(1^6)\} \ .
\end{equation}
These dimension vectors
are associated via the Kac correspondence with
minimal imaginary roots for the affine Lie algebras
$\widehat{E}_6$, $\widehat{E}_7$ and $\widehat{E}_8$,
respectively (cf. the proof of Proposition~\ref{pr:Kac indecomposables}).
(Note that the quadruple $((1^2),(1^2),(1^2),(1^2))$
that appeared in the proof of Theorem~\ref{th:k=3}
corresponds to the minimal imaginary root for $\widehat{D}_4$.)
Using formula (\ref{eq:Q-formula}),
we find $Q(\dd') = 0$ for any $\dd' \in N_{p,q,r}$.

Without loss of generality, we can assume that a triple
$\dd = (\aa, \bb, \cc)$ is reduced, i.e., all the compositions
$\aa, \bb$ and $\cc$ have all parts non-zero.
Thus, $\dd \in \Lam_{p,q,r}$, where $p, q$, and $r$ have the same meaning
as in Theorem~\ref{th:finite-type-list}.

\begin{lemma}
\label{lm:necessity}
Let $\dd \in \Lam_{p,q,r}$ be a reduced triple of compositions.
Then exactly one of the following holds:
\item {(i)} $\dd$ belongs to one of the types
$A$---$S$ in Theorem~\ref{th:finite-type-list}.
\item {(ii)} The triple $\dd$ has some $\dd' \in N_{p,q,r}$ as a summand.
\end{lemma}

\noindent {\it Proof.} Following the usual
classification of Dynkin diagrams, we present the proof in the schematic
form of a tree of implications:
\\
{\large
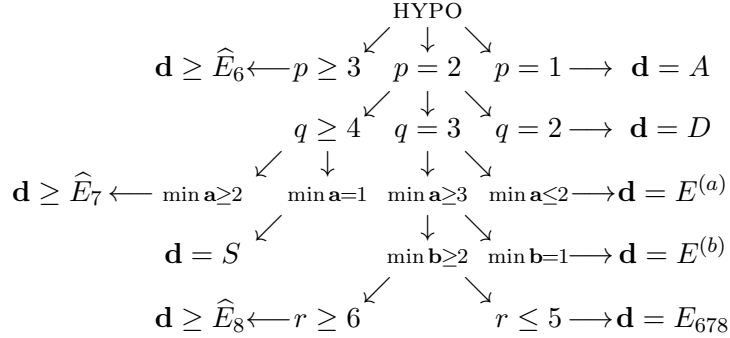
$$
\begin{array}{c@{\!}c@{\!}c@{\!}c@{\!}c@{\!}c@{\!}c@{\!}c@{\!}c@{\!}c@{\!}c}
&&&&&&\mbox{\sc hypo}&&&&\\[-.2em]
&&&&&\swarrow&\downarrow&\searrow&&&\\[-.3em]
&&\dd\geq\widehat{E}_6\,&\longleftarrow&p \geq 3&&p=2&&p=1
&\longrightarrow&\dd=A\\
&&&&&\swarrow&\downarrow&\searrow&&&\\[-.3em]
&&&&q \geq 4&&q=3&&q=2&\longrightarrow&\dd=D\\
&&&\swarrow&\downarrow&&\downarrow&\searrow&&&\\[-.3em]
\dd\geq\widehat{E}_7&\ \longleftarrow\,&
\mbox{\footnotesize $\min\aa\!\!\geq\!\! 2$}&
&\mbox{\footnotesize $\min\aa \!\!=\!\!1$}&
&\mbox{\footnotesize $\min\aa\! \!\geq\!\! 3$}&
&\mbox{\footnotesize $\min\aa\!\!\leq\!\!
2$}&\,\longrightarrow&\,\dd=E^{(a)}\\
&&&\swarrow&&&\downarrow&\searrow&&&\\[-.3em]
&&\dd=S&&&&
\mbox{\footnotesize $\min\bb \!\! \geq\!\! 2$}&
&\mbox{\footnotesize $\min\bb \!\!=\!\!1$}&\longrightarrow&\,\dd = E^{(b)}\\
&&&&&\swarrow&&\searrow&&&\\[-.3em]
&&\dd\geq\widehat{E}_8\,&\longleftarrow&r \geq 6&&&&
r \leq 5 &\longrightarrow&\,\dd = E_{678}
\end{array}
$$ }
\mbox{} \\[.5em]
The root of the tree is our
\\[.2em]
\centerline{{\sc hypothesis}: $\dd \in \Lam_{p,q,r}$
is a reduced triple of compositions, and
$1 \leq p \leq q \leq r$.}
\\[.2em]
The arrows coming from a statement
point to all possible cases resulting from the statement.
We employ the abuse of notation $\dd = A$, $\dd=D$, etc to indicate
that $\dd$ belongs to the corresponding type
in Theorem~\ref{th:finite-type-list}.
Similarly we write $\dd \geq \widehat{E}_6$, etc.~to indicate
that $\dd$ has a summand corresponding to the given affine type.
The lemma follows because every case ends in (i) or (ii),
and these conditions are clearly disjoint.
\endproof \smallskip

Combining Lemma~\ref{lm:necessity} with Proposition~\ref{pr:criterion},
we prove one direction of Theorem~\ref{th:finite-type-list}:
if $\dd$ is of finite type then it necessarily
belongs to one of the types $A$---$S$.
It remains to show that each of the conditions $A$---$S$ is
\emph{sufficient} for $(\aa, \bb, \cc)$ to be of finite type.
The following lemma is an immediate consequence of
Lemma~\ref{lm:necessity}.

\begin{lemma}
\label{lm:sufficiency-summands}
If $\dd \in \Lam_{p,q,r}$ is of one of the types $A$---$S$
then the same is true for any summand of $\dd$.
\end{lemma}

The remaining part of Theorem~\ref{th:finite-type-list} now follows
by combining Proposition~\ref{pr:criterion} with
Lemma~\ref{lm:sufficiency-summands} and the following.

\begin{lemma}
\label{lm:sufficiency-positivity}
Suppose $\dd \in \Lam_{p,q,r}$ is of one of the types $A$---$S$.
Then $Q(\dd) \geq 1$.
\end{lemma}

Thus it only remains to prove Lemma~\ref{lm:sufficiency-positivity},
which we deduce from formula
(\ref{eq:Q-formula}) and the following elementary estimates.

\begin{lemma}
\label{lm:estimates}
Let $\bb$ be a reduced composition with $|\bb| = n$ and $\ell (\bb) = q$.
Then
\item {(1)} $\Vert \bb \Vert^2 \geq n$, with equality
precisely when $\bb = (1^{n})$;
\item {(2)} if $q = 3$ then
$\Vert \bb \Vert^2 \geq 3(n-2)$, with equality
precisely when ${\rm max} (b_1, b_2, b_3) \leq 2$;
\item {(3)} if $q = 2$, and $n= 2m$ is even then
$\Vert \bb \Vert^2 \geq 2m^2$, with equality
precisely when $\bb = (m,m)$;
\item {(4)} if $q= 2$, and $n= 2m+1$ is odd then
$\Vert \bb \Vert^2 \geq 2m^2 + 2m + 1$, with equality
precisely when $\bb^+ = (m+1,m)$.
\end{lemma}

\noindent {\it Proof.} Easy. For example, part (2) is a consequence of
the identity:
$$\Vert \bb \Vert^2 -  3(n-2) =
\sum_{i=1}^3 (b_i^2 - 3b_i + 2) = \sum_{i=1}^3 (b_i - 1)(b_i - 2) \ . $$
The other parts are even simpler. \endproof \smallskip

\noindent {\it Proof of Lemma~\ref{lm:sufficiency-positivity}.}
Suppose $\dd = (\aa, \bb, \cc) \in \Lam_{p,q,r}$ is of one of the types
$A$---$S$.

\noindent {\bf Case 1.} Suppose $\dd$ is of one of the types
$A,D, E_{6}, E_7$, or $E_8$.
Then the form $Q$ is positive definite
(this is the Cartan-Killing classification), so $Q(\dd)\geq 1$.
Furthermore, the equality $Q(\dd)=1$ occurs precisely when
$\dd$ corresponds to a positive root of the associated simple Lie algebra
(cf. the proof of Proposition~\ref{pr:Kac indecomposables}).

\noindent {\bf Case 2.} Suppose $\dd$ is of type $E^{(a)}$.
Now the desired inequality $Q(\dd) \geq 1$
follows from the equality
$$\Vert \aa \Vert^2 - n^2 = 2^2 + (n-2)^2 - n^2 = 8 - 4n$$
and the inequalities $\Vert \bb \Vert^2 \geq 3(n-2)$ and
$\Vert \cc \Vert^2 \geq n$
(Lemma~\ref{lm:estimates}, parts (1), (2)).
The equality $Q(\dd)=1$ occurs precisely when
$(\aa^+,\bb^+,\cc^+) =$ $((2,2),(2,1,1),(1^4))$,\
$((3,2),(2,2,1),(1^5))$, or
$((4,2)$, $(2^3)$, $(1^6))$.

\noindent {\bf Case 3.} Suppose $\dd$ is of type $E^{(b)}$.
Let $\bb'$ be the composition obtained from $\bb$ by removing
a part equal to $1$, so that we have $|\bb'| = n-1$,
$\ell (\bb') = 2$, and
$\Vert \bb \Vert^2 = \Vert \bb' \Vert^2 + 1$.
If $n=2m$ is even then $Q(\dd) \geq 1$
follows from the inequalities
$$\Vert \aa \Vert^2 \geq 2m^2, \,\,
\Vert \bb' \Vert^2  \geq 2m^2 - 2m +1, \,\,
\Vert \cc \Vert^2 \geq 2m$$
(Lemma~\ref{lm:estimates}, parts (1), (3) and (4)).
The equality $Q(\dd)=1$ occurs precisely when
$(\aa^+,\bb^+,\cc^+) = ((m,m), (m,m-1,1),(1^{2m}))$.

If $n=2m+1$ is odd then $Q(\dd) \geq 1$
follows from the inequalities
$$\Vert \aa \Vert^2 \geq 2m^2 + 2m + 1, \,\,
\Vert \bb' \Vert^2  \geq 2m^2, \,\,
\Vert \cc \Vert^2 \geq 2m$$
(Lemma~\ref{lm:estimates}, parts (1), (3) and (4)).
The equality $Q(\dd)=1$ occurs precisely when
$(\aa^+,\bb^+,\cc^+) = ((m+1,m),(m,m,1),(1^{2m+1}))$.

\noindent {\bf Case 4.} Suppose $\dd$ is of type $S$.
Then $Q(\dd) \geq 1$
follows from
$$\Vert \aa \Vert^2 - n^2 = 1 + (n-1)^2 - n^2 = 2 - 2n$$
and $\Vert \bb \Vert^2 \geq n$,\quad
$\Vert \cc \Vert^2 \geq n$ (Lemma~\ref{lm:estimates}, part (1)).
The equality $Q(\dd)=1$ occurs precisely when
$(\aa^+,\bb^+,\cc^+) = ((n-1,1),(1^n),(1^n))$.

This completes the proofs of Lemma~\ref{lm:sufficiency-positivity}
and Theorem~\ref{th:finite-type-list}. \endproof \smallskip

As a by-product of the above argument (the examination of the
equality $Q(\dd)=1$), we immediately obtain Theorem~\ref{th:primitive}.

\section{Orbit representatives and generalized Bruhat order}
\label{sec:proofs-representatives}

\subsection{Morphisms of standard triples of flags}

We have seen (Theorem~\ref{th:indecomposables}) that the finite-type
indecomposable objects of the triple flag category $\FF_{pqr}$
occur only in dimensions $\dd \in \Pi_{pqr}$, and that
each such dimension contains a unique indecomposable isomorphism
class $I_{\dd}$.  Furthermore, by
Proposition~\ref{pr:finite indecomposables}, this $I_{\dd}$
is characterized among all objects of dimension $\dd$ by the
Schur condition $\langle I_{\dd},I_{\dd}\rangle = 1$.
Thus, to prove Theorem~\ref{th:standard form}, we need only show that
each of the standard forms listed there is Schur.
We will deduce the Schur property from a general formula for
$\langle F', F \rangle$, where $F'$ is an arbitrary object
presented in the standard form.

Let $F' = (V';A',B',C')$ be an object with $\dim V' = n$,
presented in the standard form
of Definition~\ref{def:standard form}.
Then the flags $(B',C')$ can be encoded by a family
$$\Del = ((i_1, j_1), \ldots, (i_n, j_n))$$
of $n$ pairs of indices satisfying
$1\leq i_1 \leq \cdots \leq i_n \leq q$
and $r \geq j_1 \geq \cdots \geq j_n \geq 1$.
These are defined in terms of the dimension vectors
$(b_1,\ldots,b_q)$ and $(c_1,\ldots,c_r)$ of
$B'$ and $C'$ by:
$$
i_l = \min\{ i \mid b_1 + \cdots + b_i \geq l \}
$$
$$
j_l = \min\{ j \mid c_1+\cdots+c_j \geq n\!+\!1\!-\!l \}.
$$
This means that each subspace $B'_i$ (resp. $C'_j$)
is spanned by the standard basis vectors $e_l$ such that $i_l \leq i$
(resp. $j_l \leq j$).  The set $\Del$ reflects the decomposition
of the pair of flags $B',C'$:
in the terminology of Example~\ref{ex:A-orbits},
the triple $(V';B',C') \in \FF_{q,r}$ is a direct sum
$\bigoplus_{l=1}^n I_{i_lj_l}$.

We see that the standard form
$F' = (V';A',B',C')$ is completely determined by the following
combinatorial data: a family $\Del$ and a collection
of subsets $S_1, \ldots, S_{a}$ in $\{1, \ldots, n\}$.
The sets $S_k$ must satisfy condition (3) in
Definition~\ref{def:standard form}, which says
that the set
$$K := \bigcup_{k \neq k'} (S_k \cap S_{k'}) $$
has at most two elements.

\begin{proposition}
\label{pr:standard rank}
Let $F' = (V';A',B',C')$ be an object of $\FF_{pqr}$
presented in the standard form of
Definition~\ref{def:standard form},
encoded as above by a set $\Delta$, and by subsets
$S_1, \ldots, S_{a}$.
Let $F = (V;A,B,C)$ be any object of $\FF_{pqr}$, and
define
$$
D_l = B_{i_l} \cap C_{j_l}, \qquad\qquad
E_{k} = \sum_{l \in S_k \setminus K} D_l \ .
$$
Let $v \mapsto \overline v$ denote the natural projection $V \to V/A_1$,
so that $\overline {U} = (U+A_1)/A_1 \subset V/A_1$
for any subspace $U \subset V$.
Then the dimension $\langle F', F\rangle$ of the space of
homomorphisms in $\FF_{pqr}$ from
$F'$ to $F$ is given as follows:
\smallskip

\noindent (1) If $K = \{\mu, \nu\}$
for some indices $\mu\neq\nu$, then
\begin{eqnarray}
\begin{array}{l}
\label{eq:standard rank}
\langle F', F \rangle =
\sum_{l = 1}^n \dim D_l - \dim \overline{D}_{\mu}
- \dim \overline{D}_{\nu}
- \sum_{k=1}^{a} \dim \overline{E}_{k}
\\[.1in]
\qquad\qquad
+ \dim (\overline{D}_{\mu} \cap \overline{D}_{\nu}
\cap \displaystyle\bigcap_{|S_k \cap K| = 1}
\overline{E}_{k} )
\\[.1in]
+ \dim (\displaystyle\bigcap_{K \subset S_k} \overline{E}_{k}
\ \cap \ (\, (\overline{D}_{\mu} \cap \
\displaystyle\bigcap_{S_k \cap K = \{\mu\}} \overline{E}_k)
\ + \
(\overline{D}_{\nu} \cap \
\displaystyle\bigcap_{S_k \cap K = \{\nu\}} \overline{E}_k)\, )\, ) \ .
\end{array}
\end{eqnarray}


\noindent (2) If $K = \{\mu\}$ for some index $\mu$ then
\begin{equation}
\label{eq:simple standard rank}
\langle F', F \rangle  =
\sum_{l = 1}^{n} \dim D_l - \dim \overline{D}_{\mu}
- \sum_{k=1}^{a} \dim \overline{E}_k + \dim (\overline{D}_{\mu}
\cap \hspace{-.3em}\bigcap_{\mu \in S_k}\hspace{-.3em} \overline{E}_k) \ .
\end{equation}

\noindent (3) If $K = \emptyset$ (i.e., all $S_k$ are pairwise disjoint)
then
\begin{equation}
\label{eq:very simple standard rank}
\langle F', F \rangle  =
\sum_{l = 1}^{n} \dim D_l - \sum_{k=1}^{a} \dim \overline{E}_k \ .
\end{equation}
\end{proposition}

\noindent {\it Proof.} We will only prove the most complicated
formula (\ref{eq:standard rank}).
A morphism from $F'$ to $F$ is a linear map
from $V'$ to $V$, and so is determined by the images of the basis vectors
$e_1, \ldots, e_n$; let us denote these images by $v_1, \ldots, v_n$.
By the definition, the vectors $v_l$ must satisfy the following conditions:
\begin{equation}
\label{eq:morphism conditions}
v_l \in D_l \ \ (1 \leq l \leq n) \ , \qquad
\sum_{l \in S_k} v_l \in A_1 \ \ (1 \leq k \leq a) \ .
\end{equation}
Thus $\langle F', F \rangle$ is equal
to the dimension of the subspace $U \subset V^n$ formed by
$n$-tuples $(v_1, \ldots, v_n)$ satisfying (\ref{eq:morphism conditions}).
Clearly, $U = {\rm Ker} \ (\varphi)$, where
$\varphi:\bigoplus_{l=1}^n D_l \to (V/A_1)^{a}$ is the linear map
$$\varphi_1 (v_1, \ldots, v_n) \mapsto
(\sum_{l \in S_1} \overline {v_l}, \ldots, \sum_{l \in S_{a}}
\overline {v_l})
\ .$$
Thus we have
\begin{equation}
\label{eq:rank1}
\langle F', F \rangle  =
\sum_{l=1}^n \dim D_l - {\rm rk} (\varphi) \ .
\end{equation}
Consider the subspace
$$W = \overline{D}_{\mu} \oplus \overline{D}_{\nu}
\oplus \bigoplus_{k=1}^{a} \overline{E}_k
\,\subset\, (V/A_1)^{a+2} \ .$$
Then the map $\varphi: \bigoplus_{l=1}^n D_l \to (V/A_1)^{a}$
can be factored as
$\varphi = \varphi_2 \circ \varphi_1$:
$$
\bigoplus_{l=1}^n D_l \stackrel{\varphi_1}{\longrightarrow} W
\stackrel{\varphi_2}{\longrightarrow} (V/A_1)^{a} \ ,
$$
where
$$\varphi_1: (v_1, \ldots, v_n) \mapsto
(\overline {v_\mu}, \overline {v_\nu},
\hspace{-.3em}\sum_{l \in S_1 \setminus \{\mu,\nu\}}
\hspace{-.3em} \overline {v_l}, \ldots,
\hspace{-.3em}\sum_{l \in S_{a} \setminus \{\mu,\nu\}}
\hspace{-.3em} \overline {v_l}) \ ,$$
and
$$\varphi_2: (w^{(\mu)}, w^{(\nu)}, w_1, \ldots, w_{a}) \mapsto $$
$$(\chi_1 (\mu) w^{(\mu)} + \chi_1 (\nu) w^{(\nu)} + w_1, \ldots,
\chi_{a} (\mu) w^{(\mu)} + \chi_{a} (\nu) w^{(\nu)} + w_{a}) \ .$$
(Here $\chi_k$ stands for the indicator function of the set $S_k$,
i.e., $\chi_k (l) = 1$ if $l \in S_k$, otherwise $\chi_k (l) = 0$.)
Since the sets $S_1 \setminus \{\mu, \nu\}, \ldots$,
$S_{a} \setminus \{\mu, \nu\}$ are pairwise disjoint,
the map $\varphi_1$ is surjective.
It follows that
\begin{eqnarray}
\begin{array}{l}
\label{eq:rank2}
{\rm rk} (\varphi) = {\rm rk} (\varphi_2) =
\dim W - \dim {\rm Ker} \ (\varphi_2) \\[.1in]
= \dim \overline{D}_{\mu} + \dim \overline{D}_{\nu} +
\sum_{k=1}^{a} \dim \overline{E}_k - \dim {\rm Ker} (\varphi_2) \ .
\end{array}
\end{eqnarray}

It remains to compute $\dim {\rm Ker} (\varphi_2)$.
The definition of $\varphi_2$ implies that the projection
$(w^{(\mu)}, w^{(\nu)}, w_1, \ldots, w_{a}) \mapsto (w^{(\mu)}, w^{(\nu)})$
restricts to an isomorphism between ${\rm Ker} (\varphi_2)$
and the space of pairs $(w^{(\mu)}, w^{(\nu)})$ such that
$$
w^{(\mu)} \in
\overline{D}_{\mu} \cap
\doublesubscript{\bigcap}{\mu \in S_k}{\nu \not\in S_k} \overline{E}_k,
\qquad w^{(\nu)} \in
\overline{D}_{\nu} \cap
\doublesubscript{\bigcap}{\mu \not\in S_k}{\nu \in S_k} \overline{E}_k,
\qquad w^{(\mu)} + w^{(\nu)} \in
\hspace{-.3em}\bigcap_{\{\mu,\nu\} \subset S_k}
\hspace{-.3em} \overline{E}_{k} \ .
$$
It follows that
\begin{eqnarray}
\begin{array}{l}
\label{eq:rank3}
\dim {\rm Ker} (\varphi_2) =
\dim ((\overline{D}_{\mu} \cap
\doublesubscript{\bigcap}{\mu \in S_k}{\nu \not\in S_k} \overline{E}_k)
\cap
(\overline{D}_{\nu} \cap
\doublesubscript{\bigcap}{\mu \not\in S_k}{\nu \in S_k} \overline{E}_k))
\\[.1in]
+ \dim (
\displaystyle\bigcap_{\{\mu,\nu\} \subset S_k}\hspace{-.3em} \overline{E}_{k}
\ \cap\ (
(\overline{D}_{\mu} \cap
\doublesubscript{\bigcap}{\mu \in S_k}{\nu \not\in S_k} \overline{E}_k) \ + \
(\overline{D}_{\nu} \cap
\doublesubscript{\bigcap}{\mu \not\in S_k}{\nu \in S_k} \overline{E}_k))) \ .
\end{array}
\end{eqnarray}
Combining (\ref{eq:rank1}), (\ref{eq:rank2}), and (\ref{eq:rank3})
we obtain the desired formula (\ref{eq:standard rank}).
\endproof

\subsection{Proof of Theorem~\ref{th:standard form}}
\label{sec:proof of standard form}

Let $F \in \FF_{p,q,r}$ be one of the standard-form triples of flags in
Theorem~\ref{th:standard form}.
It suffices to show that $F$ is a Schur indecomposable,
i.e., that $\langle F , F \rangle = 1$ (cf. (\ref{eq:dim Hom}) and
Proposition~\ref{pr:finite indecomposables}).

In the first and last case on the list, $\dd = ((1),(1),(1))$
and $\dd = ((1,n-1),(1^n),(1^n))$, the equality
$\langle F , F \rangle = 1$ follows at once from
(\ref{eq:very simple standard rank}).

In each of the first four $6$-dimensional cases on the list,
the desired equality $\langle F , F \rangle = 1$
is a direct consequence of (\ref{eq:standard rank}).
It is also easy to check by
an independent calculation that every morphism
from $F$ to itself is scalar.
For instance, let us do this for $\dd =((4,2),(2^3),(1^6))$.
Let $(x_{ij})$ be a $6 \times 6$ matrix that represents a morphism
$\varphi: V \to V$ in the standard basis $e_1, \ldots, e_6$.
The condition that $\varphi$ preserves the flags $B$ and $C$
means that the only non-zero matrix entries can be
$x_{11}, x_{21}, x_{22}, x_{33}, x_{43}, x_{44}, x_{55}, x_{65}$,
and $x_{66}$.
Thus we have $\varphi (e_2 + e_3) = x_{22} e_2 + x_{33} e_3 + x_{43} e_4$;
the condition that this vector lies in $A_1$ implies
that $x_{22} = x_{33}$ and $x_{43} = 0$.
Similarly, the condition that $\varphi (e_4 + e_5) \in A_1$
implies that $x_{44} = x_{55}$ and $x_{65} = 0$.
Finally, the two remaining conditions that $\varphi (e_1 + e_5)$
and $\varphi (e_2 + e_5 + e_6)$ lie in $A_1$ imply that
$x_{11} = x_{55}$, $x_{21} = 0$, and $x_{22} = x_{55} = x_{66}$.
Combining all these equalities, we see that $\varphi$ is scalar, as
desired.

For the rest of the list, the equality $\langle F , F \rangle = 1$
can be checked case by case with the help of (\ref{eq:simple standard rank}).
To simplify this procedure, we observe that all these cases satisfy
the following strengthened form of condition (3) in
Definition~\ref{def:standard form}:

\smallskip

\noindent $(3')$ Each set $S_k$ has at least two elements,
$\cup_{k=1}^{a} S_k = \{1, \ldots, n\}$,
and there exists an index $\mu$ such that
$S_k \cap S_{k'} = \{\mu\}$ for all $k\neq k'$.

\smallskip

Assuming $(3')$, we will give combinatorial conditions on subsets
$S_k$ that are necessary and sufficient for the corresponding
object $F$ to be Schur indecomposable.
This requires some terminology.

Let $F = (V;A,B,C)$ be a triple of flags with the dimension
vector $(\aa, \bb, \cc)$; let $n = \dim V$.
We associate to $\bb$ the subdivision of $[1,n] = \{1, \ldots, n\}$
into consecutive blocks $[1,b_1], [b_1 + 1, b_1 + b_2], \ldots$ of sizes
$b_1, \ldots, b_q$.
The blocks of this subdivision will be called $\bb$-blocks.
We define the $\cc$-blocks similarly, except going the
opposite way (so that the first $\cc$-block is $[n-c_1 + 1, n]$).
We say that an index $l \in [1,n]$ is \emph{$\bb$-separated}
(resp. \emph{$\cc$-separated}) from a subset $S \subset [1,n]$
if no element of $S$ smaller (resp. larger) than $l$ lies in the
same $\bb$-block (resp. $\cc$-block) with $l$.
If $l$ is both $\bb$-separated and $\cc$-separated from $S$,
we say that $l$ is \emph{$\bb \cc$-separated} from $S$.

\begin{proposition}
\label{pr:schur criterion}
Suppose $F=(V;A,B,C) \in \FF_{p,q,r}$ has the dimension
vector $(\aa, \bb, \cc)$ and is presented in a standard form
satisfying $(3')$.
Let $a=a_1$,
and denote $S'_0 = \{\mu\}$,\, $S'_k = S_k \setminus \{\mu\}$.

Then $F$ is a Schur indecomposable if and only if the subsets $S'_k$
for $k = 0, \ldots, a$ satisfy the following conditions:\\
(1) No two elements of the same $S'_k$ lie in the same $\bb$-block
or in the same $\cc$-block.\\
(2) For every two distinct indices $j$ and $k$, each $l \in S'_j$
is either $\bb$-separated or $\cc$-separated from $S'_k$.\\
(3) In the situation of (2), $S'_j$ contains an index
$\bb \cc$-separated from $S'_k$. \\
(4) For every $j$, any two elements of $S'_j$ are equivalent to each other
with respect to the equivalence relation generated by the following:
$l \sim l'$ if, for some $k \neq j$,  both $l$ and $l'$ are
$\bb \cc$-separated from $S'_k$.
\end{proposition}

\noindent {\it Proof.} We compute $\langle F , F \rangle$ using
(\ref{eq:simple standard rank}).
Note that under the condition $(3')$,
we have $E_k = \sum_{l\in S'_k} D_l$ for $k = 1,\ldots,a$.
Also define $E_0 = D_{\mu}$.  
Then formula (\ref{eq:simple standard rank}) further simplifies as follows:
\begin{equation}
\label{eq:simple standard rank2}
\langle F , F \rangle  =
\sum_{l=1}^n \dim D_l - \sum_{k=0}^a \dim \overline{E}_k
+ \dim( \bigcap_{i=0}^a \overline{E}_k)	\ .
\end{equation}
Tracing the definitions, we observe that each subspace
$D_l = B_{i_l} \cap C_{j_l}$ is spanned by all the basis vectors
$e_{l'}$ such that $l'$ is \emph{not} $\bb \cc$-separated from $\{l\}$.
In particular, $e_l \in D_l$.
It follows that $\overline{e}_{\mu}\in \overline{E}_k$ for all $k$,
and it is also clear that $\overline{e}_{\mu} \neq 0$.
Thus the last term of formula \ref{eq:simple standard rank2}
must contribute exactly 1 to $\langle F , F \rangle$,
and the first two terms must contribute 0.
We thus find that $\langle F , F \rangle=1$ if and only if:\\
(i) for each $k = 0, \ldots, a$, the sum
$\sum_{l \in S'_k} D_l$ is direct;\\
(ii) for each $k$, we have $A_1 \cap E_k = 0$; and\\
(iii) $\bigcap_{i=0}^a \overline{E_k}=\langle \overline{e}_{\mu}\rangle$, 
a one-dimensional space.

It is now completely straightforward to show that conditions (i) -- (iii)
are equivalent to conditions (1)--(4) in our proposition.
To be more precise, (i) translates into (1) and (2), (ii) translates
into (3), and (iii) into (4).
\endproof \smallskip

Now an easy inspection shows that all the remaining cases
in Theorem~\ref{th:standard form} satisfy
conditions (1)--(4) in Proposition~\ref{pr:schur criterion}.
In most of these cases, the inspection is simplified even more by the
following observation: if $\cc = (1^n)$ then condition (2) is automatic.
This completes the proof of Theorem~\ref{th:standard form}.
\endproof \smallskip

It is easy to show that in the four exceptional
$6$-dimensional cases there exist no subsets
$S_1, \ldots, S_{a}$ satisfying the conditions in
Proposition~\ref{pr:schur criterion} (this check
starts with the observation that the case $j = 0$ in condition
(3) means that an index $\mu$ must be the minimal element
of its $\bb$-block and the maximal element of its $\cc$-block).
This justifies our efforts in obtaining (\ref{eq:standard rank}).

\begin{remark}
\label{re:mutations}
{\rm  Our finite-type indecomposables $I_{\dd}$
are Schur objects of $\FF$,
also known in quiver theory as exceptional objects.
It is possible to obtain the list of representatives
in Theorem~\ref{th:standard form}
by a recursive procedure which is a special case
of the mutations of exceptional pairs studied by Rudakov, Schofield,
Crawley-Boevey, and Ringel (see \cite{ringel}).
In our situation, this procedure relies on the following
simple general proposition.
}
\end{remark}

\begin{proposition}
\label{pr:mutation}
Suppose there is a short exact sequence
$$
0\to F'
{\to} F
{\to} F'' \to 0
$$
in $\FF$ with the following properties: both $F'$ and $F''$ are
Schur indecomposables, and
$\langle F', F'' \rangle = \langle F , F' \rangle \,
\langle F'', F \rangle = 0$.
Then $F$ is a Schur indecomposable.
\end{proposition}

It turns out that, for every dimension vector $\dd$ in
Theorem~\ref{th:standard form}, one can construct a
short exact sequence as in
Proposition~\ref{pr:mutation} such that $F$ has dimension vector
$\dd$, and one of the Schur indecomposables
$F'$ and $F''$ has reduced dimension vector $((1),(1),(1))$.
The other summand
is smaller than $\dd$ and is also on our list so we can assume
that we already know its ``nice" presentation;
we can then use an explicit form of the short exact sequence
to construct a ``nice" presentation for $F$.
For instance, if $\dd = ((4,3),(3,1,3),(1^7))$ then we can choose
the dimension vectors of $F'$ and $F''$ to be respectively
$\dd' = ((3,3),(3,1,2),(1,1,0,1,1,1,1))$, and
$\dd''=((1,0),(0,0,1),(0,0,1,0,0,0,0))$.
Iterating this procedure, one can construct representatives
for all Schur indecomposables on our list
(this was in fact our original way to do it).

\subsection{Generalized Bruhat order}
\label{sec:Bruhat order}

Having determined the orbits in triple flag varieties of finite type,
we naturally ask how they fit together.
Recall that a parametrization of orbits is given by
Theorem~\ref{th:orbits} and Corollary~\ref{cor:orbit parametrization}.
We define the partial order (called \emph{degeneration order} or
\emph{generalized Bruhat order})
on the set of families $M = (m_{\dd}) \,\, (\dd \in \Pi_{p,q,r})$
by setting $M \leqdeg M'$ if the orbit $\Omega_M$
lies in the Zariski closure of $\Omega_{M'}$.

Recall from Theorem~\ref{th:indecomposables} that the orbit
$\Omega_M$ corresponds to the isomorphism class
$\bigoplus_{\dd \in \Pi_{p,q,r}} m_{\dd} I_{\dd}$
in the category $\FF_{p,q,r}$; denote this isomorphism class by $F_M$.
The following proposition is a special case of a result
due to C.~Riedtmann (cf.~\cite{riedtmann, bongartz1, bongartz2}).

\begin{proposition}
\label{pr:ranks-necessary}
If $M \leqdeg M'$ then
$\langle I_\dd , F_M \rangle \geq \langle I_\dd , F_{M'} \rangle$
for any $\dd \in \Pi_{p,q,r}$.
\end{proposition}

It would be interesting to know if the converse statement is also true,
%
i.e., if the degeneration order $M \leqdeg M'$ is given by the inequalities
$\langle I_\dd , F_M \rangle \geq \langle I_\dd , F_{M'} \rangle$ for all
$\dd \in \Pi_{p,q,r}$.
This is true when the graph $T_{p,q,r}$ is one of the Dynkin graphs
$A_n, D_n, E_6, E_7$, or $E_8$, as a consequence of general results
due to K.~Bongartz (cf.~\cite[\S4]{bongartz1}, \cite[\S5.2]{bongartz2}).

Note that Theorem~\ref{th:standard form} and formulas
(\ref{eq:standard rank}), (\ref{eq:very simple standard rank}),
and (\ref{eq:simple standard rank2}) allow us to compute
$\langle I_\dd , F_M \rangle$ explicitly for all $\dd \in \Pi_{p,q,r}$.
In particular, it is easy to compute $\langle I_\dd , I_{\dd'} \rangle$
for any two Schur indecomposables of finite type.
Knowing these numbers yields an explicit formula for
$\langle F_M, F_M \rangle$:
\begin{equation}
\label{eq:Hom with itself}
\langle F_M, F_M \rangle = \sum_{\dd, \dd' \in \Pi_{p,q,r}}
\langle I_\dd , I_{\dd'} \rangle \ m_\dd m_{\dd'} \ .
\end{equation}
This yields a formula for the (co)dimension of
any orbit $\Omega_M$.

\begin{proposition}
\label{pr:dim orbit}
The codimension of the orbit $\Omega_M$ in a multiple flag variety
${\rm Fl}_{\aa} (V) \times {\rm Fl}_{\bb} (V)  \times {\rm Fl}_{\cc} (V)$
of finite type is given by
\begin{equation}
\label{eq:codim orbit}
{\rm codim}  (\Omega_M) = \langle F_M, F_M \rangle - Q(\aa, \bb, \cc) \ .
\end{equation}
\end{proposition}

\noindent {\it Proof.} This follows at once from the formula
$$\langle F_M, F_M \rangle = \dim \mbox{Stab}_{GL(V)} (F) =
\dim GL(V) - \dim \Omega_M \ ,$$
where $F$ is any representative of $\Omega_M$
(cf. the proof of Proposition~\ref{pr:finite indecomposables}).
\endproof \smallskip

\begin{example}
\label{ex:A-orbits2}
{\rm Type $A_{q,r}$: two flags.
We use the notation of Example~\ref{ex:A-orbits}, so that
$\Omega_M$ denotes the orbit in
${\rm Fl}_{\bb} (V)  \times {\rm Fl}_{\cc} (V)$
corresponding to a $q \times r$ nonnegative
integer matrix $M = (m_{ij})$ with row sums
$b_1, \ldots, b_q$ and column sums $c_1, \ldots, c_r$.
Formula (\ref{eq:very simple standard rank})
specializes to
$$\langle I_{ij}, F \rangle = \dim (B_i \cap C_j)$$
for any pair of flags $F = (V;B,C)$.
It follows that
\begin{equation}
\label{eq:type A ranks}
\langle I_{ij}, F_M \rangle = \sum_{k=1}^i \sum_{l=1}^j m_{kl} \ .
\end{equation}
By Proposition~\ref{pr:ranks-necessary} and results of Bongartz quoted above,
the degeneration order $M \leqdeg M'$ is given by the inequalities
$$\sum_{k=1}^i \sum_{l=1}^j m_{kl} \geq
\sum_{k=1}^i \sum_{l=1}^j m'_{kl}$$
for all $i$ and $j$.
(If $\bb = \cc = (1^n)$, this is Ehresmann's original
description of the Bruhat
order on the symmetric group.)
Finally, (\ref{eq:codim orbit})
implies the following formula for the codimension of an orbit $\Omega_M$
in ${\rm Fl}_{\bb} (V)  \times {\rm Fl}_{\cc} (V)$:
\begin{equation}
\label{eq:type A codim}
{\rm codim}  (\Omega_M) = \sum_{k < i, \ l < j} m_{ij} m_{kl} \ .
\end{equation}
}
\end{example}

\begin{example}
\label{ex:S-orbits2}
{\rm Type $S_{q,r}$: two flags and a line.
We will use the notation of Example~\ref{ex:S-orbits}.
In particular, $\Omega_{\Del,M}$ denotes the orbit in
${\rm Fl}_{\aa}(V) \times {\rm Fl}_{\bb} (V)  \times {\rm Fl}_{\cc} (V)$
corresponding to a $q \times r$ nonnegative integer matrix
$M = (m_{ij})$ as above and to a non-empty set
$\Del = \{ (i_1,j_1), (i_2,j_2),$ $\ldots ,$ $(i_{t}, j_{t})\}$
such that $1\leq i_1<\ldots <i_t \leq q$, $r \geq j_1 > \ldots > j_t \geq 1$,
and $m_{ij} > 0$ for $(i,j) \in \Del$.
For any triple of flags $F = (V;A,B,C)$,
recall that $A_1$ is the only proper subspace in the flag $A$,
and $\dim A_1 = 1$.
Using (\ref{eq:very simple standard rank}), we obtain
$$\langle I_{ij}, F \rangle = \dim (B_i \cap C_j) \ ,$$
$$\langle I_{\Del'}, F \rangle =
\dim(A_1 \cap \sum_{(i,j) \in \Del'} (B_i \cap C_j)\,)
+ \sum_{(i,j)\in \Del'}\dim(B_i \cap C_j)
-\dim(\sum_{(i,j)\in \Del'}(B_i \cap C_j)\ )\ .$$
These formulas imply that
\begin{eqnarray}
\begin{array}{l}
\label{eq:type S ranks}
\langle I_{ij}, F_{\Del,M} \rangle = \sum_{k=1}^i \sum_{l=1}^j m_{kl} \ ;
\\[.1in]
\langle I_{\Del'}, F_{\Del,M} \rangle = \delta_{\Del \leq \Del'} +
\sum_{(i,j) \in {\rm In}(\Del')} \sum_{k=1}^i \sum_{l=1}^j m_{kl} \ ,
\end{array}
\end{eqnarray}
where we use the following notation: $\Del \leq \Del'$ means that
for any $(k,l) \in \Del$ there exists $(i,j) \in \Del'$
such that $k \leq i$ and $l \leq j$; the $\delta$-symbol has the usual
indicator meaning; and the operation $\Del' \mapsto {\rm In}(\Del')$
is defined by
$${\rm In}(\{ (i_1,j_1), \ldots , (i_{t}, j_{t})\}) = 
\{ (i_1,j_2), (i_2, j_3), \ldots, (i_{t-1}, j_{t})\} \ .$$
Finally, (\ref{eq:codim orbit})
implies the following formula for the codimension of an orbit
$\Omega_{\Del,M}$ in
${\rm Fl}_{\aa}(V) \times {\rm Fl}_{\bb} (V)  \times {\rm Fl}_{\cc} (V)$:
\begin{equation}
\label{eq:type S codim}
{\rm codim}  (\Omega_{\Del,M}) =
\sum_{k < i, \ l < j} m_{ij} m_{kl} +
\sum_{\{(i,j)\} \not \leq \Del} m_{ij} \ .
\end{equation}
}
\end{example}


\begin{thebibliography}{15}

\bibitem{atiyah}  M. Atiyah, On the Krull-Schmidt theorem with
application to sheaves,
{\sl Bull. Soc. Math. France} {\bf 84} (1956), 307-317.

\bibitem{bongartz1} K. Bongartz, On degenerations and extensions
of finite dimensional modules, {\sl Adv. Math.} {\bf 121} (1996), 245--287.

\bibitem{bongartz2} K. Bongartz, Degenerations for representations of
tame quivers, {\sl Ann. Sci. \'Ec. Norm. Sup.}, (4)
{\bf 28} (1995), 647--668.

\bibitem{brion} M. Brion, Groupe de Picard et nombres
caract\'eristiques des vari\'et\'es spheriques, {\sl Duke Math. J.}
{\bf 58} (1989), 397--424.

\bibitem{kac}
V. Kac, Infinite root systems, representations of graphs and invariant
theory, {\sl Invent. Math.} {\bf 56} (1980), 57--92.

\bibitem{littelmann}
P. Littelmann, On spherical double cones, {\sl J. of Algebra} {\bf 166}
(1994), 142--157.

\bibitem{riedtmann}
C. Riedtmann, Degenerations for representations
of quivers with relations, {\sl Ann. Sci. \'Ec. Norm. Sup.} (4) {\bf 19}
(1986), 275--301.

\bibitem{ringel} C. M. Ringel, Exceptional modules are tree modules,
preprint 1997.

\bibitem{simpson}  C.T. Simpson,
Products of matrices, in {\sl Differential geometry, global
analysis, and topology (Halifax, NS, 1990)}, CMS Conf. Proc. {\bf 12},
Amer. Math. Soc., Providence, RI (1991).

\end{thebibliography}
\end{document}